\documentclass[a4paper,12pt]{article}
\usepackage{mathtext}                    
\usepackage{cmap}                        
\usepackage[cp1251]{inputenc}            
\usepackage[english, russian]{babel}     
\usepackage[left=1cm,right=1cm,top=1cm,bottom=2cm]{geometry} 
\usepackage{amssymb}
\usepackage{amsmath}

\newtheorem{definition}{Definition}
\newtheorem{corollary}{Corollary}
\newtheorem{remark}{Remark}
\newtheorem{lemma}{Lemma}

\newtheorem{theorem}{Theorem}

\begin{document}

\centerline{\textbf{\LARGE{Bernstein functions of several semigroup generators}}}

\centerline{\textbf{\LARGE{ on Banach spaces under bounded perturbations. II}}}
\vspace{1cm}

\centerline{A. R. Mirotin}
\vspace{1cm}

\centerline{ Department of mathematics and programming technologies,}

 \centerline{F. Skorina Gomel State University,}

\centerline{ Gomel, 246019, Belarus}

\centerline{amirotin@yandex.ru}

\vspace{1cm}

{\small \textsc{Abstract.}  The paper deals with multidimensional Bochner-Phillips functional calculus.  In the previous paper by the author bounded perturbations of Bernstein functions of several commuting semigroup generators on Banach spaces where considered,  conditions for  Lipschitzness and estimates for the norm of commutators of such functions where proved.  Also in the one-dimensional case the Frechet differentiability of Bernstein functions of semigroup generators on Banach spaces where proved and a  generalization of Livschits-Kre\u{\i}n  trace formula derived. The aim of the present paper is to prove the Frechet differentiability of operator Bernstein functions  and the Livschits-Kre\u{\i}n  trace formula in the multidimensional setting.}

\section{Introduction}

Yu. Daletski\u{\i} and S.G. Kre\u{\i}n pioneered the study of the problem of differentiability of functions of self-adjoint operators   in  \cite{DK}. Differential calculus for functions of several commuting Hermitian operators in Hilbert spaces was studied in \cite{KPSSII}. For survey and bibliography of the theory that resulted  see the article \cite{AP1}.  It should be stressed that all these work deal with  Hilbert spaces only. The case of Banach spaces was considered in \cite{RST}, \cite{RM}, and \cite{OaM}.
In particular, the last paper  by the author was devoted to  Bernstein functions of several commuting semigroup generators on Banach spaces (they constitute the subject matter of the so called multidimensional Bochner-Phillips functional calculus).
In \cite{OaM} we gave inter alia conditions for their Lipschitzness and showed that  such functions are  $\mathcal{J}$ perturbations preserving where $\mathcal{J}$ is an arbitrary operator ideal, estimates for the norm of commutators were  also obtained; in the one-dimensional case   Frechet-differentiability and a  trace formula were  proved. The aim of the present paper is to prove the Frechet differentiability of operator Bernstein functions  and the Livschits-Kre\u{\i}n  trace formula in the multidimensional setting. (Apparently the multidimensional
Livschits-Kre\u{\i}n trace formula established in the paper is also new for
Hilbert space operators.)

 So this work could be considered as a contribution to the multidimensional Bochner-Phillips functional calculus.

As was mentioned in \cite{OaM}, the one-dimensional Bochner-Phillips functional calculus is a substantial part of the theory of
operator semigroups  and finds important applications in the theory of random
processes (see,  e.g., \cite{MirSMZ}, \cite{MirSF}). The foundations of multidimensional
calculus were laid by the author in \cite{Mir97}, \cite{Mir98}, \cite{Mir99}, \cite{SMZ2011},  \cite{IZV2015}.  Below we recall some notions and
facts from \cite{Boch}, \cite{Mir97}, \cite{Mir98}, and \cite{Mir99}, which we need for formulating our results.

\begin{definition}
 \cite{Boch} We say that a nonpositive function $\psi\in C^\infty((-\infty;0)^n)$ {\it
 belongs to the class
 ${\cal T}_n$ } (or is a nonpositive Bernstein function of $n$ variables) if all of its first partial derivatives are absolutely
monotone (a function in $C^\infty((-\infty;0)^n)$ is said to be absolutely monotone if it is nonnegative
together with its partial derivatives of all orders).
\end{definition}

Obviously,  $\psi\in{\cal T}_n$ if and only if $-\psi(-s)$ is a nonnegative Bernstein function of $n$ variables on $(0,\infty)^n$, and ${\cal T}_n$  is a cone under the pointwise addition of functions and multiplication by scalars.
As is known \cite{Boch} (see also \cite{Mir99}, \cite{mz}), each function  $\psi\in {\cal T}_n$  admits an integral representation of the form (here
and in what follows, the dot  denotes inner product in  $\mathbb{R}^n$ and the expression $s\to
-0$ means that $s_1\to -0, \ldots, s_n\to -0$)
$$
\psi(s)=c_0+c_1\cdot s+\int\limits_{\Bbb{R}_+^n\setminus \{0\}}
(e^{s\cdot v}-1)d\mu (v)\quad   (s\in(-\infty;0)^n),  \eqno(\ast)
$$
where $c_0=\psi(-0):=\lim_{s\to -0}\psi(s)$, $c_1=(c_1^j)_{j=1}^n\in \Bbb{R}_+^n$, $c_1^j=\lim_{s_j\downarrow -\infty}\psi(s)/s_j$, and $\mu$ is a positive measure on
$\Bbb{R}_+^n\setminus \{0\}$;  $\mu$
are determined by $\psi$.

A lot of  examples of  Bernstein function of one variable one can found in \cite{SSV}
(see also \cite{MirSF}, \cite{mz}).

Throughout the paper, $T_{A_1},\dots , T_{A_n}$ denote pairwise commuting one-parameter $C_0$ semigroups
(i.e., strongly continuous semigroups on $\Bbb{R}_+$) on a complex Banach space $X$ with generators $A_1, \dots ,A_n$ respectively satisfying the condition $\|T_{A_j}(t)\|\leq
M_A\quad(t\geq 0, M_A={\rm const}$)  (sometimes we write $T_{j}$ instead of $T_{A_j}$). We denote the domain of $A_j$  by $D(A_j)$ and set $A = (A_1, \dots ,A_n)$. We put $R(t, A_i):=(tI- A_i)^{-1}$ for  $t$ in $\rho(A_i)$, the resolvent set for $A_i$. Hereafter, by the commutation of operators $A_1, \dots ,A_n$
we mean the commutation of the corresponding semigroups. By ${\rm Gen}(X)$ we denote the set of all
generators of uniformly bounded $C_0$ semigroups on $X$ and by ${\rm Gen}(X)^n$, the set of all $n$-tuples  $(A_1, \dots ,A_n)$ where $A_j\in {\rm Gen}(X)$. We put also $M:=\max\{M_A, M_B\}$ for the pare $A, B$ of $n$-tuples  from  ${\rm Gen}(X)^n$. In the following $\mathcal{L}(X)$ denotes the algebra of linear bounded operators on $X$ and  $I$, the identity operator on $X$.
An operator-valued function $T_A(u) := T_{A_1}(u_1) \dots T_{A_n}(u_n) \quad (u\in
\Bbb{R}_+^n)$ is an $n$-parameter $C_0$ semigroup;
therefore, the linear manifold $D(A):=\cap_{j=1}^nD(A_j)$ is dense in $X$ \cite[ Sec. 10.10]{HiF}.

\begin{definition}
   \cite{Mir99} The value of a function  $\psi\in {\cal T}_n$ of the form ($\ast$) at $A=(A_1,\ldots ,A_n)$ applied
to $x\in D(A)$ is defined by
$$
\psi(A)x=c_0x+c_1\cdot Ax+\int\limits_{\Bbb{R}_+^n\setminus \{0\}}
(T_A(u)-I)xd\mu(u),     
$$
where $c_1\cdot Ax:=\sum_{j=1}^nc_1^jA_jx$.
\end{definition}

Given $\psi\in{\cal T}_n$ and $t\geq 0$, the function $g_t(z):=e^{t\psi(z)}$ is absolutely monotone on $(-\infty;0)^n$. It
is also obvious that $g_t(z)\leq 1$. By virtue of the multidimensional version of the Bernstein-Widder
Theorem (see, e.g.,  \cite{Boch}, \cite{BCR}), there exists a unique bounded positive measure $\nu_t$ on $\Bbb{R}_+^n$, such that, for $z\in (-\infty;0)^n$, we have
$$
g_t(z)=\int\limits_{\Bbb{R}_+^n} e^{z\cdot u}d\nu_t(u). 
$$

\begin{definition} In the notation introduced above, we set
$$
g_t(A)x=\int\limits_{\Bbb{R}_+^n} T_A(u)xd\nu_t(u)\ (x\in X) 
$$
(the integral is understood in the sense of Bochner).
\end{definition}

Obviously, $\|g_t(A)\|\leq M_A^n$.  The map $g(A):t\mapsto g_t(A)$ is a $C_0$ semigroup. In the one-dimensional case, it is called the
\textit{semigroup subordinate to} $T_A.$
In \cite{Mir09} it was noticed that the closure of the operator $\psi(A)$ exists and is the generator of the $C_0$ semigroup
$g(A)$ (cf. \cite{Mir99}.) It suggests the following final version of the definition of the operator $\psi(A)$.

\begin{definition}
\cite{Mir99} By the value of a function $\psi\in {\cal T}_n$ at an $n$-tuple $A = (A_1,\dots ,A_n)$ of
commuting operators in ${\rm Gen}(X)$ we understand the generator of the semigroup $g(A)$, i.e.,  the closure of the operator defined in the Definition 2. This value is
denoted by $\psi(A)$.
\end{definition}

The functional calculus thus arising is called multidimensional Bochner-Phillips
calculus, or $\mathcal{T}_n$-calculus.

In the sequel  we  assume for the sake of simplicity  that $c_0=c_1=0$ in the integral representation ($\ast$) of the function $\psi\in \mathcal{T}_n$ (otherwise one should replace $\psi(s)$ by $\psi(s)-c_0-c_1\cdot s$).

The notation and constraints introduced above are used in what follows without additional
explanations.

We shall use also the following results from \cite{OaM}.


\begin{theorem}
 \cite[Theorem 1]{OaM}. Let $\psi\in {\cal T}_n$. Then for every commuting families  $A=(A_1,\dots, A_n)$,  and  $B=(B_1,\dots,B_n)$ from $\mathrm{Gen}(X)^n$ such that the operators $A_i-B_i$ are bounded,  $D(A_i)=D(B_i)\ (i=1,\dots,n)$  the operator $\psi(A)-\psi(B)$ is also bounded and
 $$
 \|\psi(A)-\psi(B)\|\leq -\frac{2e}{e-1}nM^n\psi\left(-\frac{M}{2n}\|A-B\|\right),
 $$
where $\|A-B\|:=(\|A_1-B_1\|,\dots,\|A_n-B_n\|)$.
\end{theorem}

Below we shall assume that  a (two sided) operator ideal $(\mathcal{J}, \|\cdot\|_{_\mathcal{J}})$ on $X$ is \textit{symmetrically normed} in the sense that $\|ASB\|_{_\mathcal{J}}\leq \|A\|\|S\|_{_\mathcal{J}}\|B\|$
for $A, B\in\mathcal{ L}(X)$ and $S\in \mathcal{J}$.  The following theorem shows that Bernstein functions are  $\mathcal{J}$ perturbations preserving.

\begin{theorem}
\cite[Theorem 4]{OaM}. Let $(\mathcal{J}, \|\cdot\|_{_\mathcal{J}})$ be an operator ideal on $X$ and  $\psi\in {\cal T}_n$ be such that $\frac{\partial\psi}{\partial s_i}\left|_{s=-0}\right.\ne \infty$ for all  $i=1,\dots,n$. For every commuting families  $A=(A_1,\dots, A_n)$  and  $B=(B_1,\dots,B_n)$ from $\mathrm{Gen}(X)^n$ such that  $A_i-B_i\in\mathcal{ J}, D(A_i)=D(B_i)\ (i=1,\dots,n)$  the operator $\psi(A)-\psi(B)$ belongs to $\mathcal{J}$, too, and
$$
 \|\psi(A)-\psi(B)\|_{_\mathcal{J}}\leq \left. M^{n+1}\sum_{i=1}^{n}\frac{\partial\psi}{\partial s_i}\right|_{s=-0}\|A_i-B_i\|_{_\mathcal{J}}.
 $$
\end{theorem}

\section{Differentiability}
\label{2}

\begin{definition}
  (Cf. \cite{OaM}.) Let $(\mathcal{J}, \|\cdot\|_{_\mathcal{J}})$ be an operator ideal on $X$, $\psi\in {\cal T}_n$, $A$ be the $n$-tuple of pairwise commuting operators from $\mathrm{Gen}(X)$. We call the bounded linear operator   $\psi_A^\nabla:\mathcal{J}^n\to\mathcal{J}$ (transformator)   the \textit{$\mathcal{J}$-Frechet derivative} of the operator function $\psi$  at  the point $A$, if for every $n$-tuple $\Delta A\in \mathcal{J}^n$ such that $A_i+\Delta A_i\in \mathrm{Gen}(X)$ for all $i$ and operators  $A_i+\Delta A_i$  pairwise commute   we have
$$
\|\psi(A+\Delta A)-\psi(A)-\psi_A^\nabla(\Delta A)\|_{_\mathcal{J}}=o(\|\Delta A\|_{_\mathcal{J}}) \mbox{ as }  \|\Delta A\|_{_\mathcal{J}}:=\sum_{i=1}^n\|\Delta A_i\|_{_\mathcal{J}}\to 0.
$$
\end{definition}

Evidently, the Frechet derivative at    the point $A$ is unique.

Before we formulate our first result note that  if $\partial\psi/\partial s_i|_{s=-0}\ne\infty$ the derivative  $\partial\psi/\partial s_i$ of a function  $\psi\in {\cal T}_n$ equals to
$$
\int\limits_{\Bbb{R}_+^n\setminus \{0\}}e^{s\cdot v}v_id\mu(v)\quad (s\in (-\infty,0]^n)
$$
 and the measure $v_id\mu(v)$ is finite. So for every  $n$-tuple $A$ of pairwise commuting operators $A_i\in \mathrm{Gen}(X)$ the operator
$$
\frac{\partial\psi (A)}{\partial s_i}:=\int\limits_{\Bbb{R}_+^n\setminus \{0\}}T_A(v)v_id\mu(v)
$$
exists and belongs to $\mathcal{L}(X)$.

\begin{theorem}
Let $\psi\in {\cal T}_n$, and $\forall i\  \omega_i:= \partial\psi/\partial s_i|_{s=-0}\ne\infty.$ Then for every  $n$-tuple $A$ of pairwise commuting operators $A_i\in \mathrm{Gen}(X)$
the $\mathcal{ L}(X)$-Frechet derivative for the operator function $\psi$  at    the point $A$ exists and
$$
\psi_A^\nabla(C)=\sum\limits_{i=1}^n\frac{\partial\psi (A)}{\partial s_i}C_i\eqno(1)
$$
for every $n$-tuple $C\in\mathcal{ L}(X)^n$.
\end{theorem}

Proof.
The case $n=1$ was considered in \cite{OaM}. Now let $n\geq 2.$ For the proof we need the following generalization of Theorem 5 from \cite{mz} on divided differences.

\begin{lemma}
  Let the function $\psi\in {\cal T}_n$ has the integral representation ($\ast$), and $\omega_i=\partial\psi/\partial s_i|_{s=-0}\ne\infty$ \ ($i=1,\dots,n$).
Then   the function
$$
 \varphi_i(s, s_{n+1}) :=\left\{
\begin{array}
{@{\,}r@{\quad}l@{}}
\frac{\psi(s)-\psi(s_1,\dots,s_{i-1},s_{n+1},s_{i+1},\dots,s_n)}{s_i-s_{n+1}}-\omega_i, \quad {\rm if }\ s_i\ne s_{n+1},\\
\frac{\partial\psi (s)}{\partial s_i}-\omega_i, \quad {\rm if }\   s_i=s_{n+1}
\end{array}
\right.
$$
($s=(s_1,\dots,s_n)$) belongs to  ${\cal T}_{n+1}$ and has the integral representation

$$
\varphi_i(s,s_{n+1})=\int\limits_{\mathbb{R}^{n+1}_+\setminus\{0\}}
(e^{s\cdot u+s_{n+1}u_{n+1}}-1) d\mu_i(u),
$$
\textit{where $d\mu_i(u)$ is the image of the measure $1/2d\mu(v)dw$ under the mapping}
$u_j=v_j$ \textit{if} $j\ne i, n+1,$ $u_i=(v_i+w)/2,$  $u_{n+1}=(v_i-w)/2$.
\end{lemma}

Proof.
  Assume that $s_i\ne s_{n+1}.$ Then, putting $\hat{s_i}:=(s_1,\dots,s_{i-1},s_{i+1},\dots,s_n),$ we get from ($\ast$) that
$$
\frac{\psi(s)-\psi(s_1,\dots,s_{i-1},s_{n+1},s_{i+1},\dots,s_n)}{s_i-s_{n+1}}-\omega_i=
\int\limits_{\mathbb{R}^{n}_+\setminus\{0\}}
\frac{e^{s_iv_i}-e^{s_{n+1}v_i}}{s_i-s_{n+1}}e^{\hat{s_i}\cdot\hat{v_i}}d\mu(v)-\omega_i.\eqno(2)
$$
Since
$$
\frac{e^{s_iv_i}-e^{s_{n+1}v_i}}{s_i-s_{n+1}}e^{\hat{s_i}\cdot\hat{v_i}}=
\frac{1}{2}\int\limits_{-v_i}^{v_i}\left(e^{s_i\frac{v_i+w}{2}+s_{n+1}\frac{v_i-w}{2}+\hat{s_i}\cdot\hat{v_i}}-1\right)dw+v_i,
$$
formula (2) implies that
$$
\varphi_i(s,s_{n+1})=
\frac{1}{2}\int\limits_{\mathbb{R}^{n}_+\setminus\{0\}}\int\limits_{-v_i}^{v_i}
\left(e^{s_i\frac{v_i+w}{2}+s_{n+1}\frac{v_i-w}{2}+\hat{s_i}\cdot\hat{v_i}}-1\right)dwd\mu(v)
$$
$$
=\frac{1}{2}\int\limits_{\Omega_i}
\left(e^{s_i\frac{v_i+w}{2}+s_{n+1}\frac{v_i-w}{2}+\hat{s_i}\cdot\hat{v_i}}-1\right)dwd\mu(v),
$$
where $\Omega_i=\{(v,w):v\in \mathbb{R}_+^n\setminus\{0\}, w\in[-v_i,v_i]\}.$

Making the change of variables $u_j=v_j$   for $j\neq i, n+1,$ $u_i=(v_i+w)/2,$  $u_{n+1}=(v_i-w)/2$ in the last integral, we get
$$
\varphi_i(s,s_{n+1})=
\int\limits_{\mathbb{R}^{n+1}_+\setminus\{0\}}
(e^{s\cdot u+s_{n+1}u_{n+1}}-1) d\mu_i(u),\eqno(3)
$$
where $ d\mu_i(u)$ is the image of the measure $1/2d\mu(v)dw$ under the mapping
$u_j=v_j$ if $j\ne i, n+1,$ $u_i=(v_i+w)/2,$  $u_{n+1}=(v_i-w)/2.$

The case $s_i=s_{n+1}$  follows from the first assertion of the lemma as $s_{n+1}\to s_i.$

Now we claim that for every commuting $A_1,\dots, A_{n+1}\in {\rm Gen}(X)$, and for every $i$ such that $A_i-A_{n+1}\in \mathcal{L}(X)$ the following equality holds for $x\in D(A)$
$$
\varphi_i(A_1,\dots, A_{n+1})(A_i-A_{n+1})x=\psi(A_1,\dots,A_n)x-\psi(A_1,\dots,A_{i-1},A_{n+1},A_{i+1},\dots,A_n)x
$$
$$
-\left.\frac{\partial\psi}{\partial s_i}\right|_{s=-0}(A_i-A_{n+1})x.\eqno(4)
$$
For the proof first note that
 in view of Definition 2 and formula (3) for $x\in D(A)$ we have
$$
\varphi_i(A_1,\dots, A_{n+1})(A_i-A_{n+1})x=\int\limits_{\mathbb{R}^{n+1}_+\setminus\{0\}} (T_1(u_1)\dots T_{n+1}(u_{n+1})-I)(A_i-A_{n+1})xd\mu_i(u)
$$
(for simplicity we write $T_j$ instead of  $T_{A_j}$). Let  $\Omega_i$ be as in the proof of Lemma 1. If we put in the last integral  $v_j=u_j (j\ne i, n+1),$ $v_i=u_i+u_{n+1},$ $w=u_i-u_{n+1},$ then $(v,w)$ runs over $\Omega_i$ and
 we get
$$
\varphi_i(A_1,\dots, A_{n+1})(A_i-A_{n+1})x
$$
$$
= \frac{1}{2}\int\limits_{\mathbb{R}^{n}_+\setminus\{0\}}  d\mu(v)\int\limits_{-v_i}^{v_i}
\left(T_i\left(\frac{v_i+w}{2}\right)T_{n+1}\left(\frac{v_i-w}{2}\right)\prod\limits_{1\leq j\leq n, j\ne i}T_j(u_j)-I\right)(A_i-A_{n+1})xdw
$$
$$
=\int\limits_{\mathbb{R}^{n}_+\setminus\{0\}} \prod\limits_{1\leq j\leq n, j\ne i}T_j(u_j)\frac{1}{2} \int\limits_{-v_i}^{v_i}
T_i\left(\frac{v_i+w}{2}\right)T_{n+1}\left(\frac{v_i-w}{2}\right)(A_i-A_{n+1})xdwd\mu(v)
$$
$$
-\left.\frac{\partial\psi}{\partial s_i}\right|_{s=-0}(A_i-A_{n+1})x.\eqno(5)
$$
Because of the following identity  \cite[p. 211]{OaM}  ($x\in D(A)$)
$$
\frac{1}{2} \int\limits_{-v_i}^{v_i}
T_i\left(\frac{v_i+w}{2}\right)T_{n+1}\left(\frac{v_i-w}{2}\right)(A_i-A_{n+1})xdw
=(T_i(u_i)-T_{n+1}(u_i))x,
$$
 formula (5) implies (4).

Now putting  $A_{n+1}-A_i=\Delta A_i$ in the formula (3)  ($\Delta A_i\in \mathcal{L}(X)$), we have for $x\in D(A)$
$$
\psi(A_1,\dots,A_{i-1},A_{i}+\Delta A_i,A_{i+1},\dots,A_n)x-\psi(A_1,\dots,A_n)x
$$
$$
=\left.\varphi_i(A_1,\dots,A_n,A_{i}+\Delta A_i)\Delta A_ix+\frac{\partial\psi}{\partial s_i}\right|_{s=-0}\Delta A_ix. \eqno(6)
$$
By Theorem 1 the operator
$$
\alpha_i(\Delta A_i):=\varphi_i(A_1,\dots,A_n,A_{i}+\Delta A_i)-\varphi_i(A_1,\dots,A_n,A_{i})
$$
is bounded and
$$
\|\alpha_i(\Delta A_i)\|\leq -\frac{2e}{e-1}(n+1)M^{n+1}\varphi_i\left(-\frac{M}{2n+2}(0,\dots,0,\|\Delta A_i\|)\right)\to
$$
$$
 -\frac{2e}{e-1}(n+1)M^{n+1}\varphi_i(0,\dots,0)=0 \mbox{ as } \|\Delta A_i\|\to 0.
$$
Thus the formula (6)  entails the equality
$$
\psi(A_1,\dots,A_{i-1},A_{i}+\Delta A_i,A_{i+1},\dots,A_n)x-\psi(A_1,\dots,A_n)x
$$
$$
=\left.\varphi_i(A_1,\dots,A_n,A_{i})\Delta A_ix+\frac{\partial\psi}{\partial s_i}\right|_{s=-0}\Delta A_ix+\alpha_i(\Delta A_i)\Delta A_ix. \eqno(7)
$$
Let us show that
$$
\varphi_i(A_1,\dots,A_n,A_{i})x=\left.\frac{\partial\psi(A)}{\partial s_i}x-\frac{\partial\psi}{\partial s_i}\right|_{s=-0}x.\eqno(8)
$$
 To this end note that the Definition 2 implies in view of Lemma 1 that for $x\in D(A)$
$$
\varphi_i(A_1,\dots,A_n,A_{i})x
$$
$$
=\int\limits_{\mathbb{R}^{n+1}_+\setminus\{0\}}(T_1(u_1)\dots T_{i-1}(u_{i-1})
T_i(u_i+u_{n+1})T_{i+1}(u_{i+1})\dots T_{n}(u_{n})-I)x d\mu_i(u).
$$
 If we put here   $v_j=u_j (j\ne i, n+1),$ $v_i=u_i+u_{n+1},$ $w=u_i-u_{n+1},$ as in the proof of the formula (5),  we get
$$
\varphi_i(A_1,\dots,A_n,A_{i})x=\int\limits_{\Omega_i}(T_1(v_1)\dots T_{n}(v_{n})-I)x\frac{1}{2}d\mu(v)dw
$$
$$
=\frac{1}{2}\int\limits_{\mathbb{R}^{n}_+\setminus\{0\}}d\mu(v)
\int\limits_{-v_i}^{v_i}(T_1(v_1)\dots T_{n}(v_{n})-I)xdw
$$
$$
=\int\limits_{\mathbb{R}^{n}_+\setminus\{0\}}T_1(v_1)\dots T_{n}(v_{n})v_id\mu(v)-\int\limits_{\mathbb{R}^{n}_+\setminus\{0\}}v_id\mu(v)=\left.\frac{\partial\psi(A)}{\partial s_i}x-\frac{\partial\psi}{\partial s_i}\right|_{s=-0}x.
$$
This completes the proof of formula (8). Since $D(A)$ is dense in $X$,  formula (7) implies in view of (8) that ($i=1,\dots,n$)
$$
\psi(A_1,\dots,A_{i-1},A_{i}+\Delta A_i,A_{i+1},\dots,A_n)-\psi(A_1,\dots,A_n)=\frac{\partial\psi(A)}{\partial s_i}\Delta A_i+o(\|\Delta A_i\|).
$$
It follows that
$$
\psi(A+\Delta A)-\psi(A)
=\psi(A_1+\Delta A_1,\dots,A_n+\Delta A_n)-\psi(A_1,A_2+\Delta A_2,\dots,A_n+\Delta A_n)+\dots
$$
$$
+\psi(A_1,\dots,A_{n_1},A_n+\Delta A_n)-\psi(A_1,A_2,\dots,A_n)
$$
$$
=\sum\limits_{i=1}^n\frac{\partial\psi}{\partial s_i}(A_1,\dots,A_i,A_{i+1}+\Delta A_{i+1},\dots,A_n+\Delta A_n)\Delta A_i+o(\|\Delta A\|).
$$

To complete the proof it suffices  to show that
$$
\frac{\partial\psi}{\partial s_i}(A_1,\dots,A_i,A_{i+1}+\Delta A_{i+1},\dots,A_n+\Delta A_n)=
\frac{\partial\psi}{\partial s_i}(A_1,\dots,A_n)+\beta_i(\Delta A),\eqno(9)
$$
where $\|\beta_i(\Delta A)\|\to 0$ as $\|\Delta A\|\to 0\ (i=1,\dots,n).$

But
$$
\frac{\partial\psi}{\partial s_i}(A_1,\dots,A_i,A_{i+1}+\Delta A_{i+1},\dots,A_n+\Delta A_n)=
$$
$$
\int\limits_{\Bbb{R}_+^n\setminus \{0\}}T_{A_1}(v_1)\dots T_{A_i}(v_i)T_{A_{i+1}+\Delta A_{i+1}}(v_{i+1})\dots T_{A_n+\Delta A_n}(v_n)v_id\mu(v),\eqno(10)
$$
and it is known that $\|T_{A_j+\Delta A_j}(t)-T_{A_j}(t)\|\to 0$ as  $\|\Delta A_j\| \to 0$ (see, e.g., \cite[Theorem 13.5.8]{HiF}). Since all semigroups $T_{A_j}$ are bounded and  measures $v_jd\mu(v)$  are finite, formula (9) follows from
Lebesgue dominated convergence theorem  \cite[Ch. IV, Subsection 3.7, Corollary of Theorem 6]{Burb}. This finishes the proof of Theorem 3.

Note that the condition $\left.\forall i\ \partial\psi/\partial s_i\right|_{s=-0}\ne\infty$ is also necessary for the Frechet differentiability of the function $\psi$ at every point $A$ (take $A=(O,\dots,O)$) but in the case of exponentially stable semigroups  the following
corollary holds.

\begin{corollary}
Let $\psi\in {\cal T}_n$.
Then for every $n$-tuple $A$ of pairwise commuting operators from $\mathrm{Gen}(X)$ such that  $\forall i\ \|T_{A_i}(t)\|\leq Me^{\omega t}$ with
$\omega<0,$ the $\mathcal{ L}(X)$-Frechet derivative for the operator function $\psi$ at the point $A$ exists and (1) holds.

\end{corollary}

Proof.
 To use Theorem 3 we need the condition $\left.\forall i\ \partial\psi/\partial s_i\right|_{s=-0}\ne \infty$. To bypass it we apply Theorem 3 to the function $\psi(s_1+\omega,\dots,s_n+\omega)$ from $\mathcal{T}_n$ and  to the $n$-tuple $(A_1-\omega I,\dots,A_n-\omega I)$ from $\mathrm{Gen}(X)^n$.

\begin{theorem}
 Let $(\mathcal{J}, \|\cdot\|_{\mathcal{J}})$ be a symmetrically normed operator ideal on $X$,  $\psi\in {\cal T}_n$, and  $\left.\partial\psi/\partial s_i\right|_{s=-0}\ne\infty,$  $\left.\partial^2\psi/\partial s_i^2\right|_{s=-0}\ne\infty\ (i=1,\dots,n).$  For every $n$-tuple $A$ of pairwise commuting operators from $\mathrm{Gen}(X)$
the $\mathcal{J}$-Frechet derivative  for the operator function $\psi$ at the point $A$  exists and (1) holds for every $n$-tuple $C\in\mathcal{ J}^n$.
\end{theorem}

 Proof. As in the proof of previous theorem one can assume that $n\geq 2.$ We proceed as in the proof of Theorem 3 with $\Delta A_i\in \mathcal{J},$ $\|\Delta A_i\|_{\mathcal{J}}\to 0.$ Then Theorem 2 implies that the operator
$\alpha_i(\Delta A_i)$ belongs to $\mathcal{J}$, too and for $\|\Delta A_i\|_{\mathcal{J}}\to 0$ we have
$$
\left.\|\alpha_i(\Delta A_i)\|_{\mathcal{J}}\leq M^{n+1}\frac{\partial\varphi_i}{\partial s_{n+1}}\right|_{s=-0}\|\Delta A_i\|_{\mathcal{J}}
=\left.M^{n+1}\frac{1}{2}\frac{\partial^2\psi}{\partial s_i^2}\right|_{s=-0}\|\Delta A_i\|_{\mathcal{J}}\to 0.\eqno(11)
$$
Indeed,
$$
\left.\frac{\partial\varphi_i}{\partial s_{n+1}}\right|_{s=-0}:= \lim\limits_{s\to-0}\frac{\partial\varphi_i(s_1,\dots,s_{n+1})}{\partial s_{n+1}}
$$
$$
=\lim\limits_{s\to-0}\frac{(\psi(s)-\psi(s_1,\dots,s_{i-1},s_{n+1},s_{i+1},\dots,s_n))-
\frac{\partial\psi(s_1,\dots,s_{i-1},s_{n+1},s_{i+1},\dots,s_n)}{\partial s_i}(s_i-s_{n+1})}{(s_i-s_{n+1})^2}.
$$

But, by the Taylor's formula  (below $\xi$ lies between  $(s_1,\dots,s_{i-1},s_{n+1},s_{i+1},\dots,s_n)$ and $s$),
$$
\psi(s_1,\dots,s_{i-1},s_{n+1},s_{i+1},\dots,s_n)-\psi(s)=\frac{\partial\psi(s)}{\partial s_i}(s_{n+1}-s_i)+\frac{1}{2}\frac{\partial^2\psi(\xi)}{\partial s_i^2}(s_{n+1}-s_i)^2.
$$
Hence,
$$
\left.\frac{\partial\varphi_i}{\partial s_{n+1}}\right|_{s=-0}=\lim\limits_{s\to-0}\frac{\left(\frac{\partial\psi(s)}{\partial s_i}-\frac{\partial\psi(s_1,\dots,s_{i-1},s_{n+1},s_{i+1},\dots,s_n)}{\partial s_i}\right)(s_{n+1}-s_i)-\frac{1}{2}\frac{\partial^2\psi(\xi)}{\partial s_i^2}(s_{n+1}-s_i)^2}{(s_{n+1}-s_i)^2}.
$$
Applying the Taylor's formula to the first summand of the numerator in the right-hand side  we deduce from the last equality that
$$
\left.\frac{\partial\varphi_i}{\partial s_{n+1}}\right|_{s=-0}=\left.\frac{1}{2}\frac{\partial^2\psi}{\partial s_i^2}\right|_{s=-0}
$$
and the formula (11) follows.

Now  consider
$$
\beta_i(\Delta A)=\frac{\partial\psi}{\partial s_i}(A_1,\dots,A_i,A_{i+1}+\Delta A_{i+1},\dots,A_n+\Delta A_n)-
\frac{\partial\psi}{\partial s_i}(A_1,\dots,A_n)
$$
 (see formula (9)).
Since \cite[Theorem 13.4.1]{HiF}
$$
T_{A_j+\Delta A_j}(t)=\sum\limits_{m=0}^\infty S_m(t),
$$
where
$$
S_0=T_{A_j}, S_m(t)=\int\limits_{0}^t T_{A_j}(t-\tau)\Delta A_j S_{m-1}(\tau)d\tau\ (m\geq 1),
$$
we have for  $\|\Delta A_j\|_{\mathcal{J}}\to 0$
$$
\|T_{A_j+\Delta A_j}(t)-T_{A_j}(t)\|_{\mathcal{J}}\leq \left(\sum\limits_{m=1}^\infty\int\limits_{0}^t \|T_{A_j}(t-\tau)\| \|S_{m-1}(\tau)\|d\tau\right)\|\Delta A_j\|_{\mathcal{J}}\to 0
$$
(the series in the right-hand side converges, as the proof of Theorem 13.4.1 in \cite{HiF} shows). So the formula (10) implies in view of Lebesgue dominated convergence theorem that $\|\beta_j(\Delta A)\|_{\mathcal{J}}\to 0$ as $\|\Delta A_j\|_{\mathcal{J}}\to 0.$
The remaining part of the proof is the same as in  Theorem 3.

In context of Theorem 4 there is an analog of Corollary 1  for exponentially stable semigroups, as well.

\begin{corollary}
Let $(\mathcal{J}, \|\cdot\|_{_\mathcal{J}})$ be  a symmetrically normed operator ideal on $X$, $\psi\in {\cal T}_n$.
Then for every $n$-tuple $A$ of pairwise commuting operators from $\mathrm{Gen}(X)$ such that  $\forall i\ \|T_{A_i}(t)\|\leq Me^{\omega t}$ with
$\omega<0,$ the $\mathcal{J}$-Frechet derivative for the operator function $\psi$ at the point $A$ exists and (1) holds for every $n$-tuple $C\in\mathcal{ J}^n$.
\end{corollary}

The proof of this corollary is similar to the proof of corollary 1.

\section{Trace formula}
\label{3}

As is well known, the trace formula for a trace class perturbation of a self-adjoint operator  was
proved in a special case in \cite{L}  and in the general case in
 \cite{Kr1}. A survey of farther developments (in  context of  Hilbert spaces) and bibliography one can fined in \cite{BY}, \cite{Pel09}, see also recent papers \cite{Pel16}, \cite{AP2}, \cite{MN}, \cite{MNP}, and \cite{MNP2}.

 In this section we introduce a spectral shift function and prove  a Livschits-Kre\u{\i}n trace formula for a trace class perturbations of  generators of  $C_0$-semigroups on Banach space with approximation property if this semigroups are  holomorphic in the right half-plane and have a polynomial growth. Recall  that if the Banach space $X$ has the approximation property there is a continuous linear functional $\mathrm{tr}$ of norm 1  (\textit{a trace}) on the operator ideal $(\mathfrak{S}_1, \|\cdot \|_{\mathfrak{S}_1})$ of nuclear operators on $X$ (see, e. g., \cite[p. 64]{DF}). In the following $\partial^{\alpha}\psi$ denotes a multi-index derivative of a function $\psi,$
  $$\mathbb{C}_+^n:=\{z\in \mathbb{C}^n: \mathrm{Re}(z_j)>0,\ j=1,\dots,n\}.
  $$

\begin{theorem}
Let the Banach space $X$ has the approximation property.  Let  $A$ and $B$ be $n$-tuples of generators of pairwise commuting bounded  $C_0$-semigroups $T_{A_j}$ and $T_{B_j}$ respectively on $X$ holomorphic in the  half plane $\mathbb{C}_+$ and
  satisfying $\|T_{A_j}(\zeta)\|,$  $\|T_{B_j}(\zeta)\|$$\leq M|\zeta|^{m_j}$ for some $m_j\in \mathbb{Z}_+$ ($\zeta\in \mathbb{C}_+, j=1,\dots,n$). If $\forall j\  A_j-B_j\in\mathfrak{ S}_1$ there exists a unique  distribution $\eta_{A,B}$ supported in $\mathbb{R}_+^n$  such that for every $\psi\in \mathcal{T}_n$ with $\left.\partial^{2m+1}\psi\right|_{s=-0}\ne\infty$ ($m=(m_1,\dots,m_n)$) we have
$$
\mathrm{tr}(\psi(A)-\psi(B))=\int\limits_{\mathbb{R}_+^n\setminus\{0\}}\langle\eta_{A,B}(t),e^{-u\cdot t}\rangle d\mu(u),
$$
where  (as above)  $\mu$ stands for the  representing measure of $\psi$ and  $\langle\eta_{A,B}(t), e^{-u\cdot t}\rangle$ denotes  the Laplace transform of  $\eta_{A,B}.$
In particular,
$$
\mathrm{tr}(T_A(v)-T_B(v))=\langle\eta_{A,B}(t),e^{-v\cdot t}\rangle\ (v\in \mathbb{R}_+^n\setminus\{0\}).
$$

\end{theorem}

Proof.  Consider the function
$$
F(z):=T_A(z)-T_B(z)\quad (z\in \mathbb{C}_+^n).
$$

It is easy to verify that
$$
T_A(z)-T_B(z)=\left(\prod\limits_{i=1}^{n-1}T_{A_i}(z_i)\right)(T_{A_n}(z_n)-T_{B_n}(z_n))
$$
$$
+\left(\prod\limits_{i=1}^{n-2}T_{A_i}(z_i)\right)(T_{A_{n-1}}(z_{n-1})-T_{B_{n-1}}(z_{n-1}))T_{B_n}(z_n)
$$
$$
+\dots + T_{A_1}(z_1)(T_{A_2}(z_2)-T_{B_2}(z_2))\prod\limits_{i=3}^nT_{B_i}(z_i)+
(T_{A_1}(z_1)-T_{B_1}(z_1))\prod\limits_{i=2}^nT_{B_i}(z_i).\eqno(12)
$$
Theorem 2 implies that  $T_{A_i}(z_i)-T_{B_i}(z_i)\in \mathfrak{S}_1$. So,  $F: \mathbb{C}_+^n \to \mathfrak{S}_1$ by formula (12).

But  for all $x\in D(A), \mathrm{Re}(z_i)>0$
$$
(T_{A_i}(z_i)-T_{B_i}(z_i))x=\int\limits_{[0,z_i]}\frac{d}{ds}(T_{B_i}(z_i-s)T_{A_i}(s)x)ds
$$
$$
=\int\limits_{[0,z_i]}T_{B_i}(z_i-s)(A_i-B_i)T_{A_i}(s)xds. \eqno(13)
$$

 Since for $s\in [0,z_i]$
 $$
 \|T_{B_i}(z_i-s)(A_i-B_i)T_{A_i}(s)\|_{\mathfrak{S}_1} \leq M^2(|z_i-s||s|)^{m_i}\|A_i-B_i\|_{\mathfrak{S}_1}
 $$
 $$
\leq M^2(|z_i-s|+|s|)^{2m_i}\|A_i-B_i\|_{\mathfrak{S}_1}=M^2|z_i|^{2m_i}\|A_i-B_i\|_{\mathfrak{S}_1},\eqno(14)
 $$
 both sides in (13) are bounded,  formula (13) holds for all $x\in X$, and
$$
\|T_{A_i}(z_i)-T_{B_i}(z_i)\|_{\mathfrak{S}_1}\leq M^2|z_i|^{2m_i+1}\|A_i-B_i\|_{\mathfrak{S}_1}.
$$
Now it follows from (12) that
$$
\|F(z)\|_{\mathfrak{S}_1}\leq M^{n+1}\max\limits_i\|A_i-B_i\|_{\mathfrak{S}_1}\prod\limits_{i=1}^n|z_i|^{2m_i+1}.
$$
Therefore if we put
$$
f(z):=\mathrm{tr }F(z),
$$
then
$$
|f(z)|\leq \mathrm{const} \prod\limits_{i=1}^n|z_i|^{2m_i+1}. \eqno(15)
$$

We clame that $f$ is holomorphic in $\mathbb{C}_+^n$. In view of Hartogs Theorem it suffices to prove that $f$ is separately holomorphic. To simplify the notation we shall show that $f$ is holomorphic in $z_1.$ Indeed, formula (12) yields that
$$
F(z)=(T_{A_1}(z_1)-T_{B_1}(z_1))S_1(z_2,\dots,z_n)+T_{A_1}(z_1)S_2(z_2,\dots,z_n)
$$
for some  operators $S_1(z_2,\dots,z_n)\in \mathcal{L}(X)$ and $S_2(z_2,\dots,z_n)\in\mathfrak{ S}_1.$

Then for every  $z_1\in \mathbb{C}_+$ and sufficiently small $\Delta z_1$ we have
$$
F(z_1+\Delta z_1,z_2,\dots,z_n)-F(z)
$$
$$
=((T_{A_1}(z_1+\Delta z_1)-T_{B_1}(z_1+\Delta z_1))-(T_{A_1}(z_1)-T_{B_1}(z_1)))S_1(z_2,\dots,z_n)
$$
$$
+(T_{A_1}(z_1+\Delta z_1)-T_{A_1}(z_1))S_2(z_2,\dots,z_n).\eqno(16)
$$
Formula (13) implies the equality
$$
(T_{A_1}(z_1+\Delta z_1)-T_{B_1}(z_1+\Delta z_1))-(T_{A_1}(z_1)-T_{B_1}(z_1))
$$
$$
=\int\limits_{[0,z_1+\Delta z_1]}T_{B_1}(z+\Delta z_1-s)(A_1-B_1)T_{A_1}(s)ds-\int\limits_{[0,z_1]}T_{B_1}(z_1-s)(A_1-B_1)T_{A_1}(s)ds
$$
$$
=\int\limits_{[0,z_1]}T_{B_1}(z_1+\Delta z_1-s)(A_1-B_1)T_{A_1}(s)ds-\int\limits_{[0,z_1]}T_{B_1}(z_1-s)(A_1-B_1)T_{A_1}(s)ds
$$
$$
+\int\limits_{[z_1,z_1+\Delta z_1]}T_{B_1}(z_1+\Delta z_1-s)(A_1-B_1)T_{A_1}(s)ds
$$
$$
=\left(T_{B_1}\left(\frac{z_1}{2}+\Delta z_1\right)-T_{B_1}\left(\frac{z_1}{2}\right)\right)\int\limits_{[0,z_1]}T_{B_1}\left(\frac{z_1}{2}-s\right)(A_1-B_1)T_{A_1}(s)ds
$$
$$
+\int\limits_{[z_1,z_1+\Delta z_1]}T_{B_1}(z_1+\Delta z_1-s)(A_1-B_1)T_{A_1}(s)ds.
$$
Taking into account formula (14) we have
$$
\|T_{A_1}(z_1+\Delta z_1)-T_{B_1}(z_1+\Delta z_1))-(T_{A_1}(z_1)-T_{B_1}(z_1)\|_{\mathfrak{S}_1}
$$
$$
\leq \left\|T_{B_1}\left(\frac{z_1}{2}+\Delta z_1\right)-T_{B_1}\left(\frac{z_1}{2}\right)\right\|M^2\left|\frac{z_1}{2}\right|^{2m_1}\|A_1-B_1\|_{\mathfrak{S}_1}|z_1|
$$
$$
+M^2|z_1+\Delta z_1|^{2m_1}\|A_1-B_1\|_{\mathfrak{S}_1}|\Delta z_1|\to 0\quad  (\Delta z_1\to 0).
$$
Now, formula (16) shows that  the map $F: \mathbb{C}_+^n \to \mathfrak{S}_1,$ and consequently  the function $f,$ are  continuous in $z_1$ (above we used the fact that holomorphic semigroups $T_{A_1}$ and $T_{B_1}$  are norm continuous  on $\mathbb{C}_+$).

 Moreover, since $z_1\mapsto F(z_1,z_2,\dots,z_n)$ is analytic in the  half plane $\mathbb{C}_+$ with respect to the operator norm, we have for every closed path $C$ located at this half plane that
$$
\oint\limits_C f(z_1,z_2,\dots,z_n)dz_1=\mathrm{tr}\oint\limits_C F(z_1,z_2,\dots,z_n)dz_1=0.
$$
So by the Morera's Theorem the function $z_1\mapsto f(z_1,z_2,\dots,z_n)$ is  analytic in the right half plane, as well. Now, since $f$ is  analytic in $\mathbb{C}_+^n$ and satisfies (15), there is a unique  distribution $\eta_{A,B}$ supported in $\mathbb{R}_+^n$  such that $f(z)=\langle\eta_{A,B}(t),e^{-z\cdot t}\rangle$, the Laplace transform of $\eta_{A,B}$ (see, e.g., \cite[Theorem 8.13.3]{ML}).

Since, by our hypothesis,
$$
\left.\partial^{2m+1}\psi\right|_{s=-0}=\int\limits_{\mathbb{R}_+^n\setminus\{0\}}\prod\limits_{i=1}^nu_i^{2m_i+1}d\mu(u)\ne\infty
$$
and (15) holds, Definition 2 and Theorem 2
 imply that
$$
\mathrm{tr}(\psi(A)-\psi(B))=\int\limits_{\mathbb{R}_+^n\setminus\{0\}}\mathrm{tr}(T_A(u)-T_B(u))d\mu(u)
$$
$$
=\int\limits_{\mathbb{R}_+^n\setminus\{0\}}f(u)d\mu(u)=\int\limits_{\mathbb{R}_+^n\setminus\{0\}}\langle\eta_{A,B}(t),e^{-u\cdot t}\rangle d\mu(u).
$$

Putting $\mu=\delta_v,$ the Dirac measure centered at $v,$ we get the last conclusion of the theorem.

\begin{corollary}
Under the conditions of theorem 5
$$
\eta_{A,B}=L^{-1}\mathrm{tr}(T_A-T_B),
$$
where $L^{-1}$ denotes the inverse of the $n$-dimensional Laplace transform in a distributional setting.
\end{corollary}

\begin{corollary}
Under the conditions of theorem 5 for $\lambda\in \mathbb{C}_+^n$ let $R(\lambda, A):=\prod_{i=1}^n R(\lambda_i, A_i).$ Then $R(\lambda, A)-R(\lambda, B)\in \mathfrak{S}_1$ and
$$
\mathrm{tr}(R(\lambda, A)-R(\lambda, B))=\int\limits_{\mathbb{R}_+^n}\langle\eta_{A,B}(t),e^{-u\cdot t}\rangle e^{-u\cdot\lambda} du,
$$
the iterated  Laplace transformation of  $\eta_{A,B}$ ($\mathrm{the\  Stieltjes\ transform}$).
\end{corollary}

Proof.
First note that $\mathbb{C}_+\subseteq\rho(A_i).$  Since
$$
 R(\lambda_i, A_i)=\int\limits_{\mathbb{R}_+}T_{A_i}(s)e^{-\lambda_i s}ds\ (\lambda_i\in \mathbb{C}_+; i=1,\dots,n),
$$
we have
$$
R(\lambda, A)-R(\lambda, B)=\int\limits_{\mathbb{R}_+^n}(T_A(u)-T_B(u))e^{-\lambda\cdot u}du.\eqno(17)
$$

Theorem 2 with $\psi(s)=e^{s\cdot u}-1,\ \mathcal{J}=\mathfrak{S}_1$ implies that $T_A(u)-T_B(u)\in \mathfrak{S}_1$  and
$$
\|T_A(u)-T_B(u)\|_{\mathfrak{S}_1}\leq M^{n+1}\sum\limits_{i=1}^nu_i \|A_i-B_i\|_{\mathfrak{S}_1}.\eqno(18)
$$
It follows that
$$
\int\limits_{\mathbb{R}_+^n}\|T_A(u)-T_B(u)\|_{\mathfrak{S}_1}e^{-\lambda\cdot u}du<\infty.
$$
Therefore in view of (17) we get  $R(\lambda, A)-R(\lambda, B)\in \mathfrak{S}_1$ and
$$
\mathrm{tr}(R(\lambda, A)-R(\lambda, B))=\int\limits_{\mathbb{R}_+^n}\mathrm{tr}(T_A(u)-T_B(u))e^{-\lambda\cdot u}du=\int\limits_{\mathbb{R}_+^n}\langle\eta_{A,B}(t),e^{-u\cdot t}\rangle e^{-u\cdot\lambda} du
$$
by theorem 5.

\begin{remark}
 For the $n$-dimensional Stieltjes transform of distributions see \cite[Ch.10]{ML}, especially Theorem 10.8.1 therein, and \cite{Pathak}.
\end{remark}

\begin{remark}
It was shown in \cite{Pel16}
that Livschits-Kre\u{\i}n trace formula  holds for arbitrary pairs  of not necessarily bounded
self-adjoint operators with trace class difference if and only if the corresponding function is operator Lipschitz.   In \cite{MNP2}, \cite{MNP}  Livschits-Kre\u{\i}n trace formulae (for operator Lipschitz functions)  was extended to the case of pairs
of maximal dissipative ($m$-dissipative) operators  and pairs of contractions on Hilbert space. Since by \cite[Corollary 2]{OaM} (see also \cite[Corollary 13.9]{SSV}) every negative Bernstein function $\psi$ in one variable such that $\psi'(-0)\ne\infty$ is operator Lipschitz in the class of generators of contractive $C_0$-semigroups (and, more generally, in any class of generators of uniformly bounded $C_0$-semigroups with common upper bound $M$),  theorem 5 (in the case $n=1$) and \cite[Theorem 8]{OaM} are consistent with results for Hilbert space operators mentioned above.

 It should me mentioned also that when $n=1$ and $T_A$ and $T_B$ are Hilbert space contractive semigroups the result of previous corollary is closely related to \cite[Theorem 3.14]{MN}.
\end{remark}

\begin{corollary}
Under the conditions of theorem 5 let bounded one-parameter $C_0$-semigroups $g_t(A)=T_{\psi(A)}(t)$ and $g_t(B)=T_{\psi(B)}(t)$ satisfy
  $\|T_{\psi(A)}(\zeta)\|,$  $\|T_{\psi(B)}(\zeta)\|$$\leq M|\zeta|^{k}$  for some $k\in \mathbb{Z}_+$ ($\zeta\in \mathbb{C}_+$). Then
  $$
  \eta_{\psi(A),\psi(B)}=L^{-1}_{s}\int\limits_{\mathbb{R}_+^n}\langle\eta_{A,B}(t),e^{-u\cdot t}\rangle d\nu_s(u),\eqno(19)
  $$
  where $L^{-1}_{s}$ denotes the inverse of the one-dimensional Laplace transform in a distributional setting with respect to $s.$
\end{corollary}

  Proof. First note that bounded $C_0$-semigroups $T_{\psi(A)}$ and $T_{\psi(B)}$ are holomorphic in the  half plane $\mathbb{C}_+$ by \cite[Theorem 7.2]{BBD}.
Since  $\psi(A)-\psi(B)\in \mathfrak{S}_1,$ we have by corollary 3 and theorem 5  that
$$
  \eta_{\psi(A),\psi(B)}=L^{-1}_{s}\mathrm{tr}(g_s(A)-g_s(B))=L^{-1}_{s}\mathrm{tr}
    \int\limits_{\mathbb{R}_+^n}(T_A(u)-T_B(u))d\nu_s(u)
    $$
    $$
    =L^{-1}_{s}
    \int\limits_{\mathbb{R}_+^n}\mathrm{tr}(T_A(u)-T_B(u))d\nu_s(u)=
    L^{-1}_{s}\int\limits_{\mathbb{R}_+^n}\langle\eta_{A,B}(t),e^{-u\cdot t}\rangle d\nu_s(u)
 $$
 (the first integral converges in the sense of Bochner  in the  $\mathfrak{S}_1$-norm in view of (18)).

\begin{remark}
Since $g_s(z):=e^{s\psi(z)}=\int_{\mathbb{R}_+^n}e^{u\cdot z}d\nu_s(u),$  formula (19) formally can be written as
$$
\eta_{\psi(A),\psi(B)}=L^{-1}_{s}\langle\eta_{A,B}(t),e^{s\psi(-t)}\rangle.
$$
\end{remark}

\begin{corollary}
(Cf. \cite[Theorem 8]{OaM}.) Let the Banach space $X$ has the approximation property.  Let  $A$ and $B$ be generators of bounded $C_0$-semigroups $T_{A}$ and $T_{B}$ respectively on $X$ holomorphic in the  half plane $\mathbb{C}_+$ and
  satisfying $\|T_{A}(\zeta)\|,$  $\|T_{B}(\zeta)\|$$\leq M|\zeta|^{m}$  $(m\in \mathbb{Z}_+, \zeta\in \mathbb{C}_+)$. If $A-B\in\mathfrak{ S}_1$ there exists a unique  distribution $\xi_{A,B}$ supported in $\mathbb{R}_+$  such that for every $\psi\in \mathcal{T}_1$ with $\left.\psi^{(2m+1)}\right|_{s=-0}\ne\infty$  we have
$$
\mathrm{tr}(\psi(A)-\psi(B))=\int\limits_{(0,+\infty)}\langle\xi_{A,B}(t),e^{-u t}\rangle ud\mu(u).
$$
In particular,
$$
\mathrm{tr}(T_A(v)-T_B(v))=\langle\xi_{A,B}(t),e^{-v t}\rangle v.
$$
\end{corollary}

Proof.
Put  $n=1$ in theorem 5 and take for $\xi_{A,B}$ the antiderivative of $\eta_{A,B}$ supported in $\mathbb{R}_+.$

To formulate our next corollaries we need some preparations. First note that the function
$\psi_\lambda(s):=\log \lambda-\log(\lambda-s)$ ($\lambda>0$) belongs to $\mathcal{T}_1$ \cite[Example 3]{MirSF}. So, for $A\in \mathrm{Gen}(X), \lambda>0$ we can put
$$
\log(\lambda I-A):=(\log\lambda) I-\psi_\lambda(A).
$$
Note also that for $A, B\in \mathrm{Gen}(X), \lambda>0$ such that $A-B$  is nuclear the operator
$$
\log(\lambda I-B)-\log(\lambda I-A)=\psi_\lambda(A)-\psi_\lambda(B)
$$
is nuclear by theorem 2.

\begin{definition}
(Cf. \cite[formula (3.25)]{BY}). Let the Banach space $X$ has the approximation property. For $A, B\in \mathrm{Gen}(X), \lambda>0$ such that $A-B$  is nuclear define the \textit{perturbation determinant of the pair} $(A,B)$ as follows
$$
\Delta_{B/A}(\lambda)=\exp \mathrm{tr}(\log(\lambda I-B)-\log(\lambda I-A)).
$$
\end{definition}

\begin{remark}
Since $\exp \mathrm{tr}S=\det\exp S$ for nuclear $S,$ one can  define perturbation determinant of a pair $(A,B)$ by the formula
\[
\Delta_{B/A}(\lambda)=\det\exp (\log(\lambda I-B)-\log(\lambda I-A)).
\]

If, in addition,  $A$ and $B$ commute, we have for $\lambda>0$
$$
\Delta_{B/A}(\lambda)=\det((\lambda I-B)(\lambda I-A)^{-1})=\det(I+(A-B)(\lambda I-A)^{-1}).
$$
Indeed, in this case using Dyson-Phillips series (see, e.g., \cite[(13.2.4)]{HiF}) it is easy to prove that
$$
\exp(\psi_\lambda(B))\exp(\psi_\lambda(A)-\psi_\lambda(B))=\exp(\psi_\lambda(A)).
$$
where $\exp (G)$ denotes $T_G(1)$ for a generator $G$ of a $C_0$-semigroup $T_G.$
Now, putting $g_t(z)=e^{t\psi_\lambda(z)}=\lambda^t(\lambda-z)^{-t}$ in definition 3,
we  get for $t=1$ from the above equality, that
$$
\exp(\psi_\lambda(A)-\psi_\lambda(B))=(\lambda I-B)(\lambda I-A)^{-1}.
$$
\end{remark}

 Another approach to the definition of perturbation determinant of a  pair of closed  operators on Hilbert space one can fined in \cite[Section 8.1]{Ya} and \cite{MN}.

We shall use also the following notion of the Stieltjes transform of distributions \cite{Pandex}.  Let $\alpha\leq 1$ be fixed real number. The space $S_\alpha$ of test functions consists of all complex-valued functions $\varphi\in C^\infty(0,\infty)$ such that
$$
p_k(\varphi):=\sup_{t\in (0,\infty)}(1+t)^\alpha\left|\left(t\frac{d}{dt}\right)^k\varphi(t)\right|<\infty \ (k\in \mathbb{Z}_+).
$$
The topology in $S_\alpha$ in determined by the family of seminorms $(p_k)_{k\in \mathbb{Z}_+}.$
For every linear functional $f$ from the dual space  $S_\alpha'$ its Stieltjes transform $F$ is defined by the rule
$$
F(z):=\left\langle f(t),\frac{1}{t+z}\right\rangle,\quad  z\in \mathbb{C}\setminus(-\infty,0].
$$
The function $F$ is holomorphic in $\mathbb{C}\setminus(-\infty,0].$

\begin{corollary}
(Cf. \cite{Kr1}, and also \cite[formula (3.1)]{BY}).  Let $X,A,$ and $B$ be as in the above corollary.  Then
$$
\log\Delta_{B/A}(\lambda)=\left\langle\xi_{A,B}(t),\frac{1}{t+\lambda}\right\rangle\quad (\lambda>0),\eqno(20)
$$
the  Stieltjes transformation of $\xi_{A,B}$. So one can  compute $\xi_{A,B}$ via inversion theorems for the Stieltjes transform (see, e.g., \cite{ML}, \cite{Pandex}). In particular,
$$
\xi_{A,B}(t)=\lim\limits_{k\to\infty}\frac{(-t)^{k-1}}{k!(k-2)!}\frac{d^{2k-1}}{dt^{2k-1}}(t^k\log\Delta_{B/A}(t))\eqno(21)
$$
(the limit is taken in the sense of distributions).
\end{corollary}

Proof.
First note that by the proof of theorem 5 and corollary 6 the Laplace transformation of $\xi_{A,B}$ exists. So, by  definition, $\xi_{A,B}(t)=e^{p_0 t}f(t)$ for some $p_0>0$ and tempered distribution  $f$ supported in $\mathbb{R}_+.$ Since $e^{-p_0 t}S_0$ embeds in the Schwartz space $\mathcal{S}(\mathbb{R}_+),$ it follows that $\xi_{A,B}\in S'_0.$

On the other hand, $\left.\psi_\lambda^{(2m+1)}\right|_{s=-0}\ne\infty,$  and
$$
\psi_\lambda(s)=\int\limits_0^\infty(e^{su}-1)u^{-1}e^{-\lambda u}du.
$$

Thus, by corollary 6
$$
\log\Delta_{B/A}(\lambda)=\mathrm{tr}(\psi_\lambda(B)- \psi_\lambda(A))=\int\limits_0^\infty\langle\xi_{A,B}(t),e^{-u t}\rangle e^{-\lambda u}du
$$
$$
=\left\langle\xi_{A,B}(t),\int\limits_0^\infty e^{-u (t+\lambda)}du\right\rangle=\left\langle\xi_{A,B}(t),\frac{1}{t+\lambda}\right\rangle\ (\lambda>0).
$$
Formula (21) follows from the real inversion theorem for the Stieltjes transform \cite{Pandex}.

\begin{corollary}
(Cf. \cite{Kr1}, and also \cite[formula 3.7]{BY}.) Let $X,A,$ and $B$ be as in the  corollary 6. The perturbation determinant of the pair $(A,B)$ has analytic continuation to $\mathbb{C}\setminus(-\infty,0]$ and
$$
\xi_{A,B}(t)=\frac{1}{2\pi i}\lim\limits_{y\downarrow 0}\log\frac{\Delta_{B/A}(-t-iy)}{\Delta_{B/A}(-t+iy)}\quad (t>0);\eqno(22)
$$
in particular, if $\xi_{A,B}$ is real-valued,  formula (22) takes the form
$$
\xi_{A,B}(s)=\frac{1}{\pi}\lim\limits_{y\downarrow 0}\mathrm{Im}(\log(\Delta_{B/A}(-s-iy))) \quad (s>0) \eqno(23)
$$
(limits in (22) and (23) are taken in the sense of distributions).
\end{corollary}

Proof.
Since the right-hand side of the formula (20) is  holomorphic in $\mathbb{C}\setminus(-\infty,0],$ define
$$
\log\Delta_{B/A}(z):=\left\langle\xi_{A,B}(t),\frac{1}{t+z}\right\rangle,\quad z\in\mathbb{C}\setminus(-\infty,0].\eqno(24)
$$
Now we can apply the complex inversion theorem for the Stieltjes transform \cite{Pandex} and get (22). Since for the real-valued $\xi_{A,B}$ formula (24) implies
$$
\mathrm{Re}(\log(\Delta_{B/A}(-s+iy)))=\mathrm{Re}(\log(\Delta_{B/A}(-s-iy))),
$$
$$
\mathrm{Im}(\log(\Delta_{B/A}(-s+iy)))=-\mathrm{Im}(\log(\Delta_{B/A}(-s-iy))),
$$
formula (23) follows from (22).

A substantial part of properties of   perturbation determinant of pairs of  operators on Hilbert space (see, e.g.,  \cite[Section 8.1]{Ya})  is valid for $\Delta_{B/A}$. For example, the following formula holds.

\begin{corollary}
(Cf. \cite[Section 8.1, formula (4)]{Ya}.) Let $X,A,$ and $B$ be as in the  corollary 6. Then
$$
\frac{\Delta_{B/A}'(z)}{\Delta_{B/A}(z)}=\mathrm{tr}(R(z,B)-R(z,A)),\ z\in\rho(A)\cap\rho(B).
$$
\end{corollary}

Proof.
Differentiating (24) we get for $z\in\rho(A)\cap\rho(B)$ in view of corollary 4
$$
\frac{\Delta_{B/A}'(z)}{\Delta_{B/A}(z)}=\left\langle\xi_{A,B}(t),\frac{d}{dz}\frac{1}{t+z}\right\rangle=
-\left\langle\xi_{A,B}'(t),\frac{1}{t+z}\right\rangle=
$$
$$
-\left\langle\eta_{A,B}(t),\frac{1}{t+z}\right\rangle=\mathrm{tr}(R(z,B)-R(z,A)).
$$

\begin{remark}
Formula (20) implies that $\lim_{\lambda\to+\infty}\Delta_{B/A}(\lambda)=1.$ It follows also from the definition 6  and corollary 8 that $\Delta_{B/A}(z)\Delta_{C/B}(z)=\Delta_{C/A}(z)$ for   $z\in \mathbb{C}\setminus(-\infty,0],$ and operators $A, B, C\in \mathrm{Gen}(X)$ such that the pairs $(A,B)$ and $(B,C)$ satisfy all the conditions of  corollary 6.
\end{remark}

\begin{corollary}
If, in addition  to the conditions mentioned in the  corollary 6, $\xi_{A,B}$ is a measure, then
$$
\mathrm{tr}(\psi(A)-\psi(B))=\int\limits_{\mathbb{R}_+}\psi^\prime(-t)d\xi_{A,B}(t).
$$
\end{corollary}

It follows from the  corollary 6 and Tonelli's Theorem.

\vspace{2cm}

\section{Acknowledgments}

 This work was financially supported by the  Fund of Fundamental Research of Republic of Belarus. Grant number $\Phi$ 17-082.


\begin{thebibliography}{99}


\bibitem{AP1}
\textsc{A.B. Aleksandrov and V.V. Peller}, \textit{Operator Lipschitz functions}, Russian Mathematical Surveys,
\textbf{71}:4 (2016), 605--702.

\bibitem{AP2}
\textsc{A.B. Aleksandrov and V.V. Peller}, \textit{Krein’s trace formula for unitary operators and operator
Lipschitz functions}, Funct. Anal and Appl., \textbf{50}:3 (2016), 167--175.

\bibitem{BBD}
\textsc{C. Berg, K. Boyadzhiev and R. deLaubenfels}, \textit{Generation of
generators of holomorphic semigroups}, J. Austral. Math. Soc.
(Series A), {\bf 55} (1993), 246--269.



\bibitem{BCR}
\textsc{Ch. Berg, J.P.R. Christensen, P. Ressel}, \textit{Harmonic
analysis on semigroups}, Grad. Texts in Math., vol.100,
Springer-Verlag, New York-Berlin, 1984.




\bibitem{BY}
\textsc{M.S. Birman and D.R. Yafaev}, \textit{The spectral shift function. The papers of M. G. Kre\u{\i}n and their
further development}, Algebra i Analiz \textbf{4} (1992), 1--44 (Russian).
English transl.: St. Petersburg Math. J. \textbf{4} (1993), 833 -- 870.

 \bibitem{Boch}
\textsc{S. Bochner}, \textit{Harmonic analysis and the theory of probabylity}, University of California Press, Berkeley and Los Angeles, 1955.

\bibitem{Burb}
\textsc{N. Bourbaki},
 \textit{Elements de mathematique. Livre VI. Integration. 2nd ed., Ch. 1 -- 9},   Hermann, Paris, 1965 -- 1969.

\bibitem{DK}
\textsc{Yu.L. Daletskii and S.G. Kre\u{\i}n}, \textit{Integration and differentiation of functions of Hermitian operators
and application to the theory of perturbations} (Russian), Trudy Sem. Functsion. Anal., Voronezh.
Gos. Univ. 1 (1956), 81–-105.


\bibitem{DF}
\textsc{A. Defant and K. Floret}, \textit{Tensor norms and operator ideals}, North-Holland, Amsterdam, 1993.



\bibitem{HiF}
\textsc{E. Hille and R. Phillips}, \textit{Functional Analysis and Semigroups}, Amer. Math. Soc., Providence,
R.I., 1957.

\bibitem{KPSSII}
\textsc{E. Kissin, D. Potapov, V. Shulman, and F. Sukochev}, \textit{Operator smoothness in Schatten norms for functions of several variables: Lipschitz conditions, differentiability
and unbounded derivations}, Proc. London Math. Soc., \textbf{105}, 4 (2012), 661--702.


\bibitem{Kr1}
\textsc{M.G. Kre\u{\i}n}, \textit{On a trace formula in perturbation theory}, Mat. Sbornik \textbf{33} (1953), 597–-626 (Russian).


\bibitem{L}
\textsc{I.M. Lifshitz}, \textit{On a problem in perturbation theory connected with quantum statistics}, Uspekhi Mat.
Nauk \textbf{7 } (1952), 171–-180 (Russian).


\bibitem{MN}
\textsc{M. Malamud, H. Neidhart,} \textit{Trace formulas for additive and non-additive perturbations,} Adv. in Math., \textbf{274} (2015), 736--832.

\bibitem{MNP}
\textsc{M. Malamud, H. Neidhart, V. Peller,} \textit{Analytic operator Lipschitz functions in the disc and a trace formula for functions of contractions,}
 Functional Analysis and its Applications, \textbf{51}, 3 (2017), 33--55, Preprint, arXiv:1705.07225 v1 [math. FA].

\bibitem{MNP2}
\textsc{M. Malamud, H. Neidhardt, V. Peller,} \textit{A trace formula for functions of
contractions and analytic operator Lipschitz functions,}
Comptes Rendus Acad. Sci. Paris, Ser. I,  \textbf{355} (2017),  806–-811.


\bibitem{OaM}
\textsc{A. R. Mirotin}, \textit{Bernstein functions of several semigroup generators on Banach spaces under bounded perturbations}, Operators and Matrices, \textbf{11},  (2017), 199--217.



\bibitem{RM}
\textsc{A. R. Mirotin}, \textit{On some functional calculus of closed operators on Banach space. III.
Some topics in perturbation theory},  Izvestiya VUZ. Matematika, \textbf{12}, (2017), 24--34 (Russian); English transl.: Russian Math., \textbf{12}, to appear.

\bibitem{MirSMZ}
\textsc{A.R. Mirotin}, \textit{On the $\mathcal{T}$-calculus of generators for $C_0$-semigroups}, Sib. Matem. Zh., \textbf{39}, 3
(1998), 571–-582; English transl.: Sib. Math. J., \textbf{39}, 3 (1998), 493–-503.


\bibitem{MirSF}
\textsc{A. R. Mirotin},  \textit{Criteria for Analyticity of Subordinate Semigroups}, Semigroup Forum,  \textbf{78}, 2 (2009),  262--275.

\bibitem{Mir97}
\textsc{A.R. Mirotin}, \textit{Functions from the Schoenberg class $\mathcal{T}$ on the cone of dissipative elements of
a Banach algebra}, Mat. Zametki, \textbf{61}, 4 (1997), 630--633; English transl.: Math. Notes, \textbf{61}, 3–4
(1997), 524–-527.


\bibitem{Mir98}
\textsc{A.R. Mirotin}, \textit{Functions from the Schoenberg class $\mathcal{T}$ act in the cone of dissipative elements
of a Banach algebra, II}, Mat. Zametki, \textbf{64}, 3 (1998), 423–-430; English transl.: Math. Notes,
\textbf{64}, 3–-4 (1998), 364–-370.


\bibitem{Mir99}
\textsc{A.R. Mirotin}, \textit{Multidimensional $\mathcal{T}$-calculus for generators of $C_0$ semigroups}, Algebra i
Analiz, \textbf{11}, 2 (1999), 142–170; English transl.: St. Petersburg Math. J., \textbf{11}, 2 (1999), 315–-335.

\bibitem{SMZ2011}
\textsc{A.R. Mirotin}, \textit{On some properties of the multidimensional Bochner-Phillips functional calculus}, Sib. Mat. Zhurnal,  \textbf{52}, 6 (2011),  1300 -- 1312; English transl.: Siberian Mathematical Journal
 \textbf{52}, 6  (2011), pp 1032--1041.

\bibitem{IZV2015}
\textsc{A.R. Mirotin}, \textit{On joint spectra of families of unbounded operators}, Izvestiya RAN: Ser. Mat., \textbf{79}, 6
(2015), 145--170; English transl.: Izvestiya: Mathematics, \textbf{79}, 6
(2015), 1235--1259.


\bibitem{mz}
\textsc{A.R. Mirotin},  \textit{Properties of Bernstein functions of several complex variables}, Mat. Zametki, \textbf{93}, 2
(2013), 257–-265; English transl.:  Math. Notes, \textbf{93}, 2 (2013).


\bibitem{Mir09}
\textsc{A.R. Mirotin}, \textit{On multidimensional Bochner-Phillips functional calculus}, Probl. Fiz. Mat.
Tekh., \textbf{1}, 1 (2009), 63–-66 (Russian).


\bibitem{ML}
 \textsc{O.P. Misra,  J.L. Lavoine}, \textit{Transform analysis of generalized functions}, North Holland, Amsterdam, 1986.


\bibitem{Pandex}
\textsc{J. N. Pandex}, \textit{On the Stieltjes transform of generalized functions}, Proc. Camb. Phil. Soc., \textbf{71}, 1 (1972), 85--96.


\bibitem{Pathak}
\textsc{R. S. Pathak}, \textit{A distributional generalized Stieltjes transformation}, Proc. Edinburgh Math. Soc., \textbf{20}, 1 (1976), 15--22.


\bibitem{Pel09}
\textsc{V. V. Peller},
 \textit{The behavior of functions of operators under perturbations.
A glimpse at Hilbert space operators}, 287--324, Oper. Theory Adv. Appl.,
\textbf{207}, Birkhauser Verlag, Basel, 2010. Preprint, http://arxiv.org/abs/0904.1761  [math. FA].


\bibitem{Pel16}
 \textsc{V.V. Peller}, \textit{The Lifshitz-Krein trace formula and operator Lipschitz functions}, Proc. Amer. Math. Soc., \textbf{144}, (2016), 5207--5215. DOI: http://dx.doi.org/10.1090/proc/13140.

\bibitem{RST}
\textsc{J. Rozendaal, F. Sukochev,
and A. Tomskova}, \textit{Operator Lipschitz functions on Banach spaces}, Studia Mathematica, \textbf{232 }, 1  (2016), 57--92.


\bibitem{SSV}
\textsc{R. Shilling, R. Song, Z. Vondracek}, \textit{Bernstein functions.
Theory and applications},  de Greyter, Berlin-New York,  2010.


\bibitem{Ya}
\textsc{D. R. Yafaev}, \textit{Mathematical scattering theory, volume 105 of Translations
of Mathematical Monographs}, American Mathematical Society, Providence,
RI, 1992.

\end{thebibliography}
\end{document}